\documentclass[letter]{amsart}
\usepackage{amsmath,amsfonts,amsthm}
\usepackage{color}

\numberwithin{equation}{section}

\theoremstyle{plain}
\newtheorem{theorem}{Theorem}[section]

\theoremstyle{remark}

\DeclareMathOperator{\dv}{div}
\DeclareMathOperator{\tr}{tr}
\DeclareMathOperator{\id}{Id}
\DeclareMathOperator{\proj}{proj}

\newcommand{\R}{\mathbb{R}}
\newcommand{\dd}{\mathrm{d}}

\usepackage{hyperref}

\begin{document}

\title{Some basic facts on the system $\Delta u - W_u (u) = 0$}
\author{Nicholas D.\ Alikakos}
\address{Department of Mathematics\\ University of Athens\\ Panepistemiopolis\\ 15784 Athens\\ Greece}
\email{\href{mailto:nalikako@math.uoa.gr}{\texttt{nalikako@math.uoa.gr}}} 

\begin{abstract}
We rewrite the system $\Delta u - W_u (u) = 0$, for $u: \R^n \to \R^n$, in the form $\dv T = 0$, where $T$ is an appropriate stress-energy tensor, and derive certain {\em a priori} consequences on the solutions. In particular, we point out some differences between two paradigms: the phase-transition system, with target a finite set of points, and the Ginzburg--Landau system, with target a connected manifold.
\end{abstract}

\maketitle

\section{Introduction}
This paper is concerned with entire solutions to the system
\begin{equation}\label{system}
\Delta u - W_u(u) = 0,\text{ for } u: \R^n \to \R^n,
\end{equation}
where $W \in C^2(\R^n; \R)$ and $W_u := ( \partial W/ \partial u_1, \ldots, \partial W / \partial u_n )^{\top}$. Two distinguished examples are: (a) the phase-transition model or vector Allen--Cahn equation, where $W$ has a finite number of global minima $a_1, a_2, \ldots , a_N$ (see Baldo \cite{baldo}, Bronsard and Reitich \cite{bronsard-reitich}) and (b) the Ginzburg--Landau system $\Delta u - (|u|^2 -1)u = 0$ (see Bethuel, Brezis, and H\'elein \cite{bethuel-brezis-helein}).

First, we introduce a stress-energy tensor $T$, that is, an $n \times n$ matrix $T=T(u,\nabla u)$, and show that \eqref{system} can be written in the form
\begin{equation}
\dv T = 0,
\end{equation}
for solutions $u \in W^{1,2}_{\mathrm{loc}} \cap L^{\infty}_{\mathrm{loc}}$.  Next, by following the general procedure due to Schoen \cite{schoen} for divergence-free tensors, we establish, under the hypothesis 
\begin{equation}\label{positivity}
W \geq 0,
\end{equation}
the (weak) monotonicity formula
\begin{equation}\label{monotonicity_formula}
\frac{\dd}{\dd R} \left( \frac{1}{R^{n-2}} \int_{|x-x_0|<R} \left( \frac{1}{2}|\nabla u|^2 + W(u) \right) \dd x \right) \geq 0,
\end{equation}
where $\nabla u$ is the matrix $( \partial u_i / \partial x_j )$, for $i,j=1, \ldots, n$, $|\cdot|$ is the Euclidean norm, $R>0$, and $x_0$ is an arbitrary point in $\R^n$. Finally, as an application, we obtain a Liouville-type result: For a solution $u \in W^{1,2}_{\mathrm{loc}} \cap L^{\infty}_{\mathrm{loc}}$ to system \eqref{system} under the hypothesis \eqref{positivity} and for $n \geq 2$, we show that
\begin{equation}
\int_{\R^n} \left( \frac{1}{2} |\nabla u|^2 + W(u) \right) \dd x < \infty \quad \text{implies} \quad u \equiv \mathrm{constant}.
\end{equation}

Formula \eqref{monotonicity_formula} corresponds to the usual monotonicity formula for harmonic maps (see Schoen and Uhlenbeck \cite{schoen-uhlenbeck}). The (strong) monotonicity formula 
\begin{equation}\label{strong_monotonicity_formula}
\frac{\dd}{\dd R} \left( \frac{1}{R^{n-1}} \int_{|x-x_0|<R} \left( \frac{1}{2}|\nabla u|^2 + W(u) \right) \dd x \right) \geq 0
\end{equation}
implies \eqref{monotonicity_formula}. It can be motivated by the relationship between the solutions of the phase-transition model and minimal surfaces, and thus to the monotonicity of the quantity 
\[
\frac{1}{R^{n-1}}\int_{|x-x_0|<R} |\nabla u| \, \dd x.
\]
However, the validity of \eqref{strong_monotonicity_formula} depends on the estimate 
\begin{equation}
|\nabla u|^2 \leq 2 W(u).
\end{equation}
This last estimate is a well-known result of Modica \cite{modica1} for solutions to the scalar Allen--Cahn equation
\begin{equation}
\Delta u - W'(u) = 0,\text{ for } u: \R^n \to \R,
\end{equation}
but it is not known\footnote{See Section \ref{remarks:section} for useful information on this point, which we owe to Alberto Farina.} for systems like \eqref{system}. In the case of graphs $u:\R^n \to \R$, Modica \cite{modica2} has also established \eqref{strong_monotonicity_formula}. Caffarelli, Garofalo, and Segala \cite{caffarelli-garofalo-segala} have extended both results of Modica to a wider class of scalar equations under the hypothesis that $W \geq 0$, with $W=0$ on a finite set of points. Formula \eqref{monotonicity_formula} generalizes analogous results of Bethuel, Brezis, and Orlandi \cite{bethuel-brezis-orlandi}, Rivi\`ere \cite{riviere}, and Farina \cite{farina2} for the Ginzburg--Landau system.

Ilmanen \cite{ilmanen} has introduced a predecessor of the stress-energy tensor in his work on the evolution scalar Allen--Cahn equation versus motion by mean curvature. In the context of the Ginzburg--Landau system, we note the book of Sandier and Serfaty \cite{sandier-serfaty}, where this tensor has been utilized in a number of ways. Our source is Alikakos and Betel\'u \cite{alikakos-betelu} where the tensor was introduced specifically for writing the system as a divergence-free condition in order to derive the Plateau angle conditions.

In Section \ref{stress:section} we give a derivation of the well-known monotonicity formula for minimal surfaces (see Simon \cite{simon}, Ecker \cite{ecker}) under hypotheses of smoothness, based on the stress-energy tensor, by applying Schoen's procedure \cite{schoen}.

\section{The weak monotonicity formula}
\subsection{The stress-energy tensor}
We begin by introducing a {\em stress-energy tensor} for vector fields $u: \R^n \to \R^n$. Let 
\begin{equation}\label{tensor}
T := \frac{1}{2} 
\left( \begin{array}{c}
|u_{x_1}|^2 - \displaystyle{\sum_{i\neq 1}^{n} |u_{x_i}|^2 - 2 W(u)},~~ 2 u_{x_1} \!\cdot u_{x_2},~ \cdots~,~~ 2 u_{x_1} \!\cdot u_{x_n} \\
2 u_{x_2} \!\cdot u_{x_1},~~ |u_{x_2}|^2 - \displaystyle{\sum_{i\neq 2}^{n} |u_{x_i}|^2 - 2 W(u)},~ \cdots~,~~ 2 u_{x_2} \!\cdot u_{x_n}\\
\ddots \\
2 u_{x_n} \!\cdot u_{x_1},~~ 2 u_{x_n} \!\cdot u_{x_2},~ \cdots~,~~ |u_{x_n}|^2 - \displaystyle{\sum_{i\neq n}^{n} |u_{x_i}|^2 - 2 W(u)} \end{array} \right).
\end{equation}
So, $T=(T_1, T_2, \ldots, T_n)^\top$ and $\dv T = (\dv T_1, \dv T_2, \ldots, \dv T_n)^\top$.

The relevance of $T$ to the system \eqref{system} can be seen by calculating its divergence. We have
\begin{align*}
\dv T_j & = \sum_{k \neq j} \big(u_{x_j} \!\cdot u_{x_k}\big)_{\!,x_k} + \frac{1}{2} \bigg( |u_{x_j}|^2 -\sum_{i \neq j} |u_{x_i}|^2 - 2 W(u) \bigg)_{\!\!,x_j}\\
& = \sum_{k \neq j} \big(u_{x_j x_k} \!\cdot u_{x_k} + u_{x_j} \!\cdot u_{x_k x_k}\big) + u_{x_j} \!\cdot u_{x_j x_j} - \sum_{i \neq j} u_{x_i} \!\cdot u_{x_i x_j} - W_u \!\cdot u_{x_j} \\
& = u_{x_j} \!\cdot \big(\Delta u - W_u(u)\big),
\end{align*}
thus, 
\begin{equation}
\dv T = \left( \begin{array}{c}
u_{x_1} \!\cdot \big(\Delta u - W_u(u)\big) \smallskip\\ 
u_{x_2} \!\cdot \big(\Delta u - W_u(u)\big) \\
\vdots \\ 
u_{x_n} \!\cdot \big(\Delta u - W_u(u)\big)
\end{array} \right) = (\nabla u)^\top \big( \Delta u - W_u(u) \big).
\end{equation}

Therefore, in the class of smooth vector fields, equation \eqref{system} is almost equivalent to $T$ being divergence-free since
\begin{equation}\label{div-free}
\left\{
\begin{array}{l}
\Delta u - W_u(u) = 0 \quad \Rightarrow \quad \dv T = 0,\medskip\\
\Delta u - W_u(u) = 0 \quad \Leftarrow \quad \dv T = 0,\  \det \nabla u \neq 0.
\end{array}\right.
\end{equation}

Next, we calculate the trace of $T$.
\begin{align}\label{trace}
\tr T = \sum_{i=1}^{n} T_{ii} & = \frac{1}{2} \Bigg( \sum_{i=1}^{n} |u_{x_i}|^2 - 2nW(u) - \sum_{i=1}^{n} \bigg( \sum_{k \neq i} |u_{x_k}|^2 \bigg) \Bigg) \nonumber\\
& = \frac{1}{2} \Bigg( \sum_{i=1}^{n} |u_{x_i}|^2 - 2nW(u) - \sum_{i=1}^{n} \bigg( \sum_{k =1}^{n} |u_{x_k}|^2 - |u_{x_i}|^2 \bigg) \Bigg) \nonumber\\
& = \frac{2-n}{2} \sum_{i=1}^{n} |u_{x_i}|^2 - nW(u).
\end{align}

Finally, we introduce the {\em interface-energy density}
\begin{equation}
e(u) := \frac{1}{2} |\nabla u|^2 + W(u),
\end{equation}
where $\nabla u$ is the matrix $( \partial u_i / \partial x_j )$, for $i,j=1, \ldots, n$ and $|\cdot|$ is the Euclidean norm.

Notice that 
\begin{align}
\tr T & = -n\, e(u) + |\nabla u|^2,\label{tr1}\\
\tr T & \leq - (n-2)\, e(u).\label{tr2}
\end{align}

\subsection{The monotonicity formula}\label{monotonicity:section}
\begin{theorem}\label{th1}
Let $u$ be a $W^{1,2}_{\mathrm{loc}}(\R^n; \R^n) \cap L^{\infty}_{\mathrm{loc}}(\R^n; \R^n)$ solution to the system
\[
\Delta u - W_u (u) = 0,\text{ for } u: \R^n \to \R^n,
\]
under the hypothesis that $W$ is $C^2(\R^n; \R)$ and $W \geq 0$. Then, we have 
\begin{equation}\label{mono1}
\frac{\dd}{\dd R}\left( R^{-(n-2)} E_{B_R} (u) \right) \geq 0,\text{ for } R>0,
\end{equation}
where
\begin{equation}\label{e-def}
E_{B_R} (u) = \int_{B_R} \left( \frac{1}{2} |\nabla u|^2 + W(u) \right) \dd x,
\end{equation}
with $x_0 \in \R^n$ arbitrary and $B_R := B(x_0; R)$ the ball in $\R^n$.
\end{theorem}
\begin{proof}
We begin by noting the simple fact that
\begin{equation}\label{simple}
T + e(u) \id = (\nabla u)^\top (\nabla u) \geq 0,
\end{equation}
where $\id$ stands for the identity matrix on $\R^n$.

Take now $x_0 = 0$. Following Schoen \cite{schoen}, we have
\begin{align}\label{one}
\sum_{i,j} \int_{B_R} (x^i\, T_{ij})_{,j} &= \sum_{i,j} \int_{B_R} (\delta_{ij}\, T_{ij} + x^i\, T_{ij,j}) \nonumber\\
& = \sum_i \int_{B_R} T_{ii} \quad (\text{from } \eqref{div-free}) \nonumber\\
& \leq - (n-2) \int_{B_R} e(u) \quad (\text{from } \eqref{tr2}).
\end{align}
On the other hand, by the divergence theorem and for $\nu = x/R$,
\begin{align}\label{two}
\sum_{i,j} \int_{B_R} (x^i\, T_{ij})_{,j} & = \sum_{i,j} \int_{\partial B_R} x^i\, T_{ij}\, \nu_j \nonumber\\
& = R\, \sum_{i,j} \int_{\partial B_R} T_{ij}\, \nu_i\, \nu_j \nonumber\\
& = R \int_{\partial B_R} \langle T \nu, \nu \rangle \nonumber\\
& \geq - R \int_{\partial B_R} e(u) \quad (\text{from } \eqref{simple}) \nonumber\\
& = - R\, \frac{\dd E_{B_R} (u)}{\dd R}.
\end{align}
where $\langle \cdot\, , \cdot \rangle$ is the Euclidean inner product.

Combining \eqref{one} and \eqref{two} we obtain 
\[
-(n-2)\, E_{B_R} (u) \geq - R\, \frac{\dd E_{B_R} (u)}{\dd R},
\]
or, equivalently,
\begin{equation}\label{conclusion}
\frac{\dd}{\dd R} \left( R^{-(n-2)} E_{B_R} (u) \right) \geq 0.
\end{equation}
By utilizing the translation invariance of \eqref{system}, we conclude that $x_0$ can be arbitrary in \eqref{conclusion} and the proof of the theorem is complete.
\end{proof}

\subsection{Remarks}\label{remarks:section}
Note that if the analog of the Modica gradient estimate 
\begin{equation}\label{modica-gradient-estimate}
|\nabla u|^2 \leq 2 W(u)
\end{equation}
is assumed to hold for system \eqref{system}, then \eqref{tr2} can be strengthened to 
\begin{equation}\label{strenghened}
\tr T \leq -(n-1)\, e(u),
\end{equation}
from which the strong monotonicity formula
\begin{equation}\label{strong}
\frac{\dd}{\dd R} \left( R^{-(n-1)} E_{B_R} (u) \right) \geq 0, \text{ for } R>0,
\end{equation}
follows (by using in \eqref{one} formula \eqref{strenghened} instead of \eqref{tr2}). Formula \eqref{strong} has been derived already by Modica \cite{modica1} (see also Caffarelli, Garofalo, and Segala \cite{caffarelli-garofalo-segala}) for entire solutions to the scalar Allen--Cahn equation
\begin{equation}
\Delta u - W'(u) = 0, \text{ for } u:\R^n \to \R,\ W:\R \to \R.
\end{equation}
Theorem \ref{th1} is more general from the point of view of $W$. The estimate \eqref{modica-gradient-estimate} does not hold in such generality. Actually neither \eqref{modica-gradient-estimate} nor \eqref{strong} hold for the Ginzburg--Landau system \cite{farina2}.

The stress-energy tensor $T$ was introduced in joint unpublished work with Betel\'u \cite{alikakos-betelu} (see also the M.Sc.\ thesis of Dimitroula \cite{dimitroula} for further elaboration) where it is utilized for deriving the Plateau angle conditions and was motivated by the work of Bronsard and Reitich \cite{bronsard-reitich}. See also the related work of Sandier and Serfaty \cite{sandier-serfaty} on the Ginzburg--Landau system. The predecessor of $T$ for scalar fields $u: \R^n \to \R$ was introduced before by Ilmanen in \cite{ilmanen}. 

The physical meaning of $T$ can be seen via the scaled Allen--Cahn system
\begin{equation}\label{scaled-system}
\varepsilon \Delta u - \frac{1}{\varepsilon} W_u (u) = 0, \text{ for } u: \R^n \to \R^n.
\end{equation}
This requires a digression which we now take. As is well-known from the theory of $\Gamma$-convergence (see De Giorgi and Franzoni \cite{degiorgi-franzoni}, Modica and Mortola \cite{modica-mortola}, Modica \cite{modica3}, Sternberg \cite{sternberg1}, Baldo \cite{baldo}, and for expository work Alberti \cite{alberti} and Braides \cite{braides}),
\begin{equation}\label{gamma-limit}
\underset{\varepsilon \to 0}{\Gamma \text{-lim}} \int_{\Omega} \bigg( \frac{\varepsilon}{2} |\nabla u|^2 + \frac{1}{\varepsilon} W(u) \bigg) \dd x = F(u),\text{ for } \Omega \text{ bounded, open in } \R^n,
\end{equation}
with
\begin{equation}\label{perimeter}
F(u) = \sum_{i<j} \sigma_{ij} \mathcal{H}^{n-1} (S_{ij} u),
\end{equation}
where $S_{ij}u$ is the interface in $\Omega$ between the phases $\{ u=a_i \}$ and $\{ u=a_j \}$, $\mathcal{H}^{n-1}$ is the $(n-1)$-dimensional Hausdorff measure, and $\sigma_{ij}$ is the corresponding surface-tension coefficient (there holds, see below, $\sigma_{ij} + \sigma_{jk} \geq \sigma_{ik}$). For the vector case see Baldo \cite{baldo}, where $W \geq 0$ and vanishing precisely at $a_1, \ldots , a_N$. 

Without a constraint, the global minimizers $u^\varepsilon$ are trivial. $\Gamma$-convergence methods treat the integral constraint (see \cite{baldo}) and lead to surfaces of constant mean curvature. Local minimizers $u^\varepsilon$ have been treated in Kohn and Sternberg \cite{kohn-sternberg} and, from a point of view closer to the present article, in Hutchinson and Tonegawa \cite{hutchinson-tonegawa}. The blow-down limit in \eqref{gamma-limit} relates global minimizers $u^\varepsilon$ of 
\[
\int_\Omega \left( \frac{\varepsilon}{2} |\nabla u|^2 + \frac{1}{\varepsilon} W(u) \right) \dd x
\]
to global minimizers of the perimeter functional \eqref{perimeter}, that is, to minimal surfaces $S_{ij}$. As is well-known, for minimal surfaces the following monotonicity formula holds (see Simon \cite{simon}, Ecker \cite{ecker})
\begin{equation}\label{minimal-monotonicity}
\frac{\dd}{\dd R} \left( R^{-(n-1)}\, \mathcal{H}^{n-1} \big( S_{ij} \cap B^n (p; R) \big) \right) \geq 0, \text{ for } 0 \leq R < \delta,
\end{equation}
where $B^n(p;R)$ is any ball in $\R^n$ with center $p$, $\delta$-away from the boundary $\partial \Omega$. Notice the relationship between \eqref{strong} and \eqref{minimal-monotonicity}.

\subsection{The limiting stress-energy tensor and the monotonicity formula for minimal surfaces}\label{stress:section}
After this background, we are ready to derive formally the limiting stress-energy tensor, which has a very clear physical meaning. Let $U_{ij}$ be the connection between $a_i$ and $a_j$ (see Sternberg \cite{sternberg2}, Alikakos and Fusco \cite{alikakos-fusco-1}), that is, let
\begin{equation}\label{connection}
\left\{
\begin{array}{l}
\ddot{U}_{ij} (\eta) - W_u \big(U_{ij} (\eta)\big) = 0, \text{ where } U_{ij}: \R \to \R^n,\ \dot{U}:= \dd U / \dd \eta,\medskip\\
U_{ij} (-\infty) = a_i,\ U_{ij} (+\infty) = a_j, \text{ for } i \neq j,\  U_{ij} (0) = 0.
\end{array}\right.
\end{equation}
It is well-known that $U_{ij}$ satisfies the equipartition relation 
\[
\frac{1}{2} |\dot{U}_{ij}|^2 = W(U_{ij})
\]
and also
\[
\sigma_{ij} = \int_\R \bigg( \frac{1}{2}|\dot{U}_{ij}|^2 + W(U_{ij}) \bigg) \dd \eta.
\]
For small $\varepsilon >0$, the minimizer $u^\varepsilon$ can be approximated by
\[
u^\varepsilon \approx U \bigg( \frac{d(x,S)}{\varepsilon} \bigg),
\]
where $d(\cdot, S)$ is the signed distance from the $(n-1)$-manifold $S$, $S=S_{ij}$ is the part of the interface $I = \cup_{i<j} S_{ij}$ which realizes the distance of $x$ from $I$, and $U=U_{ij}$ is the corresponding connection. We note that this approximation is appropriate away from the intersections of $S_{ij}$'s (for example, away from triple junctions in the plane). However, the contribution to the energy of the intersection set is proportional to the size of small neighborhoods surrounding the set, and so can be ignored. (This is supported by formal results in Bronsard and Reitich \cite{bronsard-reitich} and rigorous work in Baldo \cite{baldo}, plus a lot of accumulated evidence on diffuse-interface models.)

The scaled stress-energy tensor is given by
\begin{equation}\label{scaled-tensor}
T_\varepsilon := \frac{1}{2} 
\left( \begin{array}{c}
\varepsilon|u_{x_1}|^2 - \displaystyle{\varepsilon \sum_{i\neq 1}^{n} |u_{x_i}|^2 - 2 \frac{W(u)}{\varepsilon}},~~ 2 \varepsilon u_{x_1} \!\cdot u_{x_2},~ \cdots~,~~ 2\varepsilon u_{x_1} \!\cdot u_{x_n} \\
2\varepsilon u_{x_2} \!\cdot u_{x_1},~~ \varepsilon |u_{x_2}|^2 - \displaystyle{\varepsilon \sum_{i\neq 2}^{n} |u_{x_i}|^2 - 2 \frac{W(u)}{\varepsilon}},~ \cdots~,~~ 2\varepsilon u_{x_2} \!\cdot u_{x_n}\\
\ddots \\
2\varepsilon u_{x_n} \!\cdot u_{x_1},~~ 2\varepsilon u_{x_n} \!\cdot u_{x_2},~ \cdots~,~~ \varepsilon |u_{x_n}|^2 - \displaystyle{\varepsilon \sum_{i\neq n}^{n} |u_{x_i}|^2 - 2 \frac{W(u)}{\varepsilon}} \end{array} \right),
\end{equation}
and corresponds to \eqref{scaled-system}. By evaluating $T_\varepsilon$ on $u^\varepsilon$, by making use of the approximate relationship above, and by utilizing the canonical coordinates around $S$ (see \cite[Appendix]{gilbarg-trudinger}), we obtain 
\begin{equation}
\lim_{\varepsilon \to 0} T_\varepsilon (u^\varepsilon) = T_0 := \sigma \mathcal{H}^{n-1} \llcorner S 
\left( \begin{array}{cccc}
d^{2}_{x_1} -1 & d_{x_1} d_{x_2} & \cdots & d_{x_1} d_{x_n} \medskip\\
d_{x_2} d_{x_1} & d^{2}_{x_2} -1& \cdots & d_{x_2} d_{x_n} \\
& & \ddots \\
d_{x_n} d_{x_1} & d_{x_n} d_{x_2} & \cdots & d^{2}_{x_n} -1 
\end{array} \right),
\end{equation}
in the sense of Radon measures. For example, for the first entry in the matrix above, we have
\[
\lim_{\varepsilon \to 0} \int_{\R^n} \biggl( \varepsilon |u_{x_1}|^2 - \varepsilon \sum_{i \neq 1} |u_{x_i}|^2 - 2 \frac{W(u^\varepsilon)}{\varepsilon} \biggr) \phi(x)\, \dd x = 2 \sigma (d_{x_1}^{2} -1) \int_S \phi(y)\, \dd S_y,
\]
for every $\phi \in C_{\mathrm{c}}(\R^n; \R)$.
Thus, 
\begin{equation}\label{projection}
T_0 \big|_S = \sigma ( \nabla d \otimes \nabla d - \id ),
\end{equation}
that is, $T_0$ is an orthogonal projection in $\R^n$, on the tangent space of $S$. Therefore, $T_0$ is a symmetric matrix generating tangential forces and, as such, it is properly called a stress-energy tensor since the surface tension force is supposed to be tangential to the interface (see also Ilmanen \cite{ilmanen}).

We write \eqref{projection} by dropping the subscript for simplicity,
\[
T \big|_S = \sigma ( \nabla d \otimes \nabla d - \id ),\quad T_{ij} \big|_S = \sigma (d_{x_i} d_{x_j} - \delta_{ij}).
\]

Next we derive \eqref{minimal-monotonicity} under the hypothesis of smoothness.

\subsubsection*{A}\label{step-a}
Since $T = - \sigma \proj_{T_x S}$, we have that 
\[
\dv_S T = - \sigma \mathbf{H},
\]
where $\mathbf{H} := H \nu$ is the mean-curvature vector, with $\nu$ a unit normal. Thus, $\dv_S T = 0$ is equivalent to $H=0$.

\subsubsection*{B}\label{step-b}
Utilizing that $T$ is a projection gives
\[
\dv_S T(X) = - \sigma \mathbf{H} \cdot X - \sigma \dv_S X,
\]
for any vector field $X$ in $\R^n$. Using that $H=0$, we have
\[
\int_{S \cap B_R} \dv_S T(X) = -\sigma \int_{S \cap B_R} \dv_S X,
\]
and by selecting $X(x) = x - p$ and $B_R = B_R(p)$, for $p \in S$, we obtain
\[
-\sigma \int_{S \cap B_R} \dv_S X = -\sigma (n-1)\, \mathcal{H}^{n-1} (S \cap B_R).
\]

\subsubsection*{C}\label{step-c}
By the divergence theorem on surfaces (see Simon \cite{simon}, Ecker \cite[Appendix]{ecker}),
\[
\int_{S \cap B_R} \dv_S X = - \int_{S \cap B_R} H \cdot X + \int_{S \cap \partial B_R} X \cdot \gamma,\]
where $\gamma$ a unit co-normal. Hence,
\begin{align*}
\int_{S \cap B_R} \dv_S T(X) & = - \sigma \int_{S \cap B_R} \dv_S X \\
& = -\sigma \int_{S \cap \partial B_R} (x-p) \cdot \gamma \\
& \geq  - \sigma R \,\mathcal{H}^{n-2} (S \cap \partial B_R) \\
& = - \sigma R \frac{\dd}{\dd R} \Big( \mathcal{H}^{n-1} (S \cap B_R) \Big).
\end{align*}

\subsubsection*{D}\label{step-d}
Combining {\em B.} and {\em C.} we obtain
\[
- \sigma (n-1)\, \mathcal{H}^{n-1} (S \cap B_R) \geq -\sigma R \frac{\mathrm{d}}{\mathrm{d} R} \Big( \mathcal{H}^{n-1} ( S \cap B_R ) \Big) ,
\]
or, equivalently,
\[
\frac{\mathrm{d}}{\mathrm{d} R} \Big( R^{-(n-1)}\, \mathcal{H}^{n-1} ( S \cap B_R ) \Big) \geq 0.
\]

The above gives an answer, under hypotheses of smoothness, to a question raised in Schoen \cite{schoen} on the possibility of such an approach for minimal surfaces, that is, the existence of a physically meaningful stress-energy tensor. It is this question that motivated the present work.

The application of the divergence theorem on a vector field tangent to the surface as a means for deriving the monotonicity formula is well known (see Simon \cite[p.~83]{simon}). The point in Schoen \cite{schoen} is that the stress-energy tensor, a physically motivated quantity, treats in a unified way a variety of setups. The calculation above can be extended to varifolds as in Simon \cite[p.~235]{simon}.

\section{A Liouville theorem}
\begin{theorem}
Let $u$ be a $W^{1,2}_{\mathrm{loc}}(\R^n; \R^n) \cap L^{\infty}_{\mathrm{loc}}(\R^n; \R^n)$ solution to the system
\[
\Delta u - W_u (u) = 0,\text{ for } u: \R^n \to \R^n,
\]
under the hypothesis that $W$ is $C^2(\R^n; \R)$ and $W \geq 0$. Then, for $n \geq 2$ we have that
\begin{equation}
\int_{\R^n} \left( \frac{1}{2} |\nabla u|^2 + W(u) \right) \dd x < \infty \quad \text{implies} \quad u \equiv \mathrm{constant}.
\end{equation}
\end{theorem}
\begin{proof}
For $n\geq 3$, by \eqref{mono1}, for $\sigma \geq R$, 
\begin{align*}
\frac{1}{R^{n-2}} \int_{|x| < R} \left( \frac{1}{2} |\nabla u|^2 + W(u) \right) \dd x & \leq \frac{1}{\sigma^{n-2}} \int_{|x| < \sigma} \left( \frac{1}{2} |\nabla u|^2 + W(u) \right) \dd x \\
& \leq \frac{1}{\sigma^{n-2}} \int_{\R^n} \left( \frac{1}{2} |\nabla u|^2 + W(u) \right) \dd x.
\end{align*}
By taking the limit $\sigma \to \infty$ we are set.

For $n=2$, \eqref{trace} gives
\[
\tr T = -2 W(u).
\]
Following the proof of Theorem \ref{th1}, we obtain
\[
R \frac{\dd E_{B_R} (u)}{\dd R} \geq 2 \int_{B_R} W(u) \,\dd x \geq 2 \int_{B_{R_0}} W(u) \,\dd x,
\]
for $R \geq R_0$. Hence, integrating we obtain
\[
E_{B_R} (u) \geq E_{B_{R_0}} (u) + 2 \log \frac{R}{R_0} \int_{B_{R_0}} W(u) \,\dd x.
\]
Thus, the hypothesis implies that $W(u) \equiv 0$. But then the components $u_1$, $u_2$ of $u$ are harmonic, with
\[
\int_{\R^2} | \nabla u_i |^2 \,\dd x < \infty, \text{ for } i=1,2.
\]
Therefore, $u_i \equiv \text{constant}$ by the mean-value property of harmonic functions applied to $\partial u_i / \partial x_j$. The proof of the theorem is complete.
\end{proof}

Note that an examination of the arguments above renders the estimates
\begin{equation*}
\left\{
\begin{array}{l}
E_{B_R} (u) = \mathrm{o} (R^{n-2}), \text{ as } R \to \infty,\ n \geq 3\quad \text{implies} \quad u \equiv \mathrm{constant},\medskip\\
E_{B_R} (u) = \mathrm{o} (\log R), \text{ as } R \to \infty,\ n = 2\quad \text{implies} \quad u \equiv \mathrm{constant}.
\end{array}\right.
\end{equation*}

For related work on scalar equations see Modica \cite{modica2} and Caffarelli, Garofalo, and Segala \cite{caffarelli-garofalo-segala}. For Liouville theorems for the related Ginzburg--Landau systems see Farina \cite{farina}, \cite{farina2}.

\section*{Acknowledgments}
This work was carried out during a visit at Stanford University in the spring of 2009. Thanks are due to Rafe Mazzeo and Rick Schoen for several stimulating conversations and to the Department of Mathematics of Stanford University as a whole for the great hospitality. Thanks also go to Peter Bates, Alex Freire, and Panagiotis Smyrnelis for their comments on a previous version of the paper.

\nocite{*}
\bibliographystyle{plain}

\end{document}